\documentclass{article}
\usepackage{amsmath}
\usepackage[left=2.0cm,right=2.0cm,top=2.5cm,bottom=2.5cm]{geometry}

\numberwithin{equation}{section}
\newtheorem{thm}{Theorem}[section]
\newtheorem{cor}[thm]{Corollary}
\newtheorem{prop}[thm]{Proposition}

\newtheorem{lem}[thm]{Lemma}

\newtheorem{ex}{Example}[section]

\newcommand{\be}{\begin{equation}}
\newcommand{\ee}{\end{equation}}
\newcommand{\ben}{\begin{enumerate}}
\newcommand{\een}{\end{enumerate}}

\newcommand{\qed}{\hspace*{\fill}Q.E.D.}  
\title{On Landsberg general $(\alpha,\beta)$-metrics with a conformal 1-form}
\author{Shasha Zhou and Benling Li\footnote{Corresponding author. Research is supported by the NNSFC(11371209), ZPNSFC(LY13A010013) and K.C. Wong Magna Fund in Ningbo University. }\\
{\small \it Department of Mathematics, Ningbo University }\\
{\small \it Ningbo, Zhejiang Province 315211, P.R. China }\\
{\small \it sz18211656622@126.com;  libenling@nbu.edu.cn}
}

\date{April 8, 2017}

\begin{document}

\maketitle

\begin{abstract}
\bigskip
   In this paper, we study almost regular Landsberg  general $(\alpha,\beta)$-metrics in Finsler geometry.
   The corresponding equivalent equations are given. By solving the equations, we give the classification of Landsberg general $(\alpha,\beta)$-metrics
   under the conditon that $\beta$ is closed and conformal to $\alpha$. Under this condition, we prove that regular Landsberg general $(\alpha,\beta)$-metrics must be Berwaldian when the dimension is greater than two. For the almost regular case, the classification also is
   given and some new non-Berwaldian Landsberg metrics are found.
{\bf MR(2010) Subject Classification: 53B40, 53C60}
\end{abstract}

\section{Introduction}
Finsler geometry is just Riemannian geometry without the quadratic restriction \cite{Chern}. More colorful than the Riemannian case, there are
some non-Riemannian quantities in Finsler geometry, such as the Berwald curvature, the Landsberg curvature, Douglas curvature, S-curvature and etc. They all vanish for Riemannian metrics, hence they are said to be non-Riemannian\cite{A1}.  A Finsler metric  is called a {\it Berwald metric} if its Berwald curvature vanishes and  is called a {\it Landsberg metric} if its Landsberg curvature vanishes.
On Berwald manifold (manifold with Berwald metric), all the tangent spaces  with the induced norm are linearly
isometric to each other. The local structures of Berwald metrics have been studied by Z.I. Szab\'{o} \cite{A2}. As we known that all Berwald metrics are Landsberg metrics. This is because that on Landsberg manifold (manifold with Landsberg metric), the tangent spaces with the
induced norm  are all isometric.
It is a long existing open problem in Finsler geometry:

 {\it Whether or not any Landsberg metric is a Berwald metric ?  }

Matsumoto called it the most important unsolved problem in Finsler geometry in 2004.
Due to the many unsuccessful attempts made for finding non-Berwald Landsberg metrics, D. Bao called them "unicorns".
In order to answer this problem, it is natural to study and construct non-Riemannian Landsberg metrics.
Some geometers considered this problem and found this problem is still open when the metric is regular \cite{B3, C1, V1, A3, A4}.
Thus, there is no regular Landsberg metric found which is non-Berwaldian till now.
However, in the almost regular case when the metric is allowed to be singular in some directions,  the question has a satisfied answer.
In 2006, G. Asanov proved that his metrics arising from Finslerian General Relativity are actually Landsberg metrics but not Berwald metrics
\cite{GA1, GA2}. After that, Z. Shen obtained a two-parameter
family of Landsberg metrics including Asanov's examples\cite{Z1}. Then it is meaningful to ask how to find more Landsberg metrics which
are non-Berwaldian even in the almost regular case?

The metrics obtained by G. Asanov and Z. Shen are belong to {\it $(\alpha,\beta)$-metrics} which are defined
by the following form
\[F=\alpha \phi(s),\ \ \ s=\frac{\beta}{\alpha},\]
where $\phi(s)$ is a $C^{\infty}$ function, $\alpha=\sqrt{a_{ij}(x)y^iy^j}$ is a Riemannian metric and $\beta=b_{i}(x)y^{i}$ is a $1$-form.
Randers metric is the simplest non-Riemannian Finsler metric in the form $F=\alpha+\beta$ \cite{R1}. It is well-known that a Randers metric is a Landsberg metric if and only if it is a Berwald metric \cite{Ma1, R2}. In 2009, Z. Shen proved that a regular $(\alpha,\beta)$-metric is a Landsberg metric if and only if it is a Berwald metric and constructed a family of almost regular non-Berwald Landsberg $(\alpha,\beta)$-metrics \cite{Z1}.

To study the "unicorns" problem and find more Landsberg metrics, its natural to  consider a more general metric class.
 {\it General $(\alpha,\beta)$-metric} was first introduced by C. Yu and H. Zhu in \cite{Y1}. By definition, a general $(\alpha,\beta)$-metric is a Finsler metric expressed in the following form,
$$F=\alpha \phi(b^2,s),\ \ \ s=\frac{\beta}{\alpha},$$
where $\phi(s)$ is a $C^{\infty}$ function, $\alpha=\sqrt{a_{ij}(x)y^iy^j}$ is a Riemannian metric and $\beta=b_{i}(x)y^{i}$ is a $1$-form, $b:=\|\beta_x\|_\alpha$.
It is easy to see that $F$ is the {\it spherically symmetric metric} when $\alpha=|y|$ is the Euclidean metric  and $\beta =\langle x, y\rangle$ is the Euclidean inner production. Obviously, $\beta$ is closed and conformal to $\alpha$ for all spherically symmetric metrics, i.e.,
\be \label{conformal}
b_{i|j}=c a_{ij},
\ee
where $c=c(x)$ is a scalar function. When $c=0$, then $\beta$ is said to be parallel to $\alpha$.
Some papers studied general $(\alpha,\beta)$-metrics under this condition \cite{Huang, LS1, Z1, Y1, Zhou}.
Thus, it is interesting to study the problems when $\beta$ is closed and conformal to $\alpha$.

In this paper, we mainly consider the "unicorns" problem.
We find the equivalent equations of the Berwald and Landsberg general $(\alpha,\beta)$-metrics under the condition (\ref{conformal}). Based on these equivalent equations,  a more refined characterization is obtained for the Landsberg general $(\alpha,\beta)$-metrics.
For Berwald metrics, we get the following theorem.

\begin{thm}\label{thm1.3}
Let $F=\alpha \phi(b^2,\frac{\beta}{\alpha})$ be a general $(\alpha,\beta)$-metric on an n-dimensional manifold with $n\geq 3$. Suppose that $\beta$ is closed and conformal to $\alpha$, i.e.,
\be b_{i|j}=c a_{ij},\label{bij3}
\ee
where $c=c(x)\neq0$ is a scalar function.
If $F$ is a Berwald metric, then it must be a Riemannian metric.
\end{thm}

When $c=c(x)=0$, it is obvious that $F$ is a Berwald metric with the same geodesic coefficients of $\alpha$.
In \cite{H1}, H. Zhu studied the general $(\alpha,\beta)$-metrics with isotropic Berwald curvature. She obtained the equivalent equations and find the  solutions. In fact, the solutions are just Riemannian metrics. Thus the above Theorem \ref{thm1.3} is based on Zhu's solutions.

In \cite{M1}, M. Zohrehvand and H. Maleki studied the regular case when $\beta$ is closed and conformal to $\alpha$ and found that Landsberg general $(\alpha,\beta)$-metrics must be Berwaldian in this case.
However, they didn't consider the almost regular case.
A function $F=\alpha \phi(b^2,\frac{\beta}{\alpha})$ is called an {\it almost regular general $(\alpha,\beta)$-metric} if $\beta$ satisfies that $\|\beta\|_\alpha \leq b_0$, $\forall x\in M$. If general $(\alpha,\beta)$-metrics are allowed to be singular in two extremal directions, then there are some non-trivial solutions. An almost regular general $(\alpha,\beta)$-metric $F=\alpha \phi(b^2,\frac{\beta}{\alpha})$ might be singular in the two extremal direction $y\in T_{x}M$ with $\beta(x,y)=\pm b_{0}\alpha(x,y)$.
In Proposition \label{thm1.1}, we give the equivalent equations of almost regular Landsberg general $(\alpha,\beta)$-metrics (non-Riemannian) (\ref{ZYFC1}) and (\ref{ZYFC2}).
To determine the Landsberg metrics, the efficient way is to solve these equations. Actually, we obtain the following result about the general solutions of (\ref{ZYFC1}) and (\ref{ZYFC2}).

\begin{thm}\label{thm1.2}
Let $F=\alpha \phi(b^2,\frac{\beta}{\alpha})$ be an almost regular non-Riemannian general $(\alpha,\beta)$-metric on an n-dimensional manifold with $n\geq3$. Suppose that $\beta$ is closed and conformal to $\alpha$, i.e.,
\be b_{i|j}=c a_{ij},\label{bij2}
\ee
where $c=c(x)\neq0$ is a scalar function. Then $F$ is a Landsberg metric if and only if the function $\phi=\phi(b^2,\frac{\beta}{\alpha})$ is given by
\be
\phi=\lambda_{3}e^{\int^s_{0}{\frac{2c_{1}\sqrt{b^2-t^2}+b^2t\lambda_{1}+2b^2t\lambda_{0}}{2tc_{1}\sqrt{b^2-t^2}
+b^2t^2\lambda_{1}+b^2-2b^2(b^2-t^2)\lambda_{0}}}dt}.\label{TJ}
\ee
where $c_{1}=c_{1}(b^2)$, $\lambda_{0}=\lambda_{0}(b^2)$, $\lambda_{1}=\lambda_{1}(b^2)$ and $\lambda_{3}=\lambda_{3}(b^2)>0$ are $C^{\infty}$ functions of $b^2$ satisfy the following condition,
\begin{equation}
\begin{split}
\int^s_{0}{A_{1}(b^2,t)}dt=&\frac{1}{2b^4s\lambda_{3}}\Big\{\Big[b^4(b^2-s^2)(2\lambda_{0}+s^2\lambda_{2})
-2sc_{1}(b^2-s^2)^{\frac{3}{2}}-b^4\Big]\lambda_{3}A(b^2,s)\\&+\Big[b^4s(\lambda_{2}s^2+2\lambda_{1}+2\lambda_{0})
+2c_{1}(b^2-s^2)^{\frac{3}{2}}\Big]\lambda_{3}-2b^4s\lambda_{3}'\Big\},\label{TJTJ}
\end{split}
\end{equation}
where
\be
A(b^2,t)=\frac{2c_{1}\sqrt{b^2-t^2}+b^2t\lambda_{1}+2b^2t\lambda_{0}}{2tc_{1}\sqrt{b^2-t^2}+b^2t^2\lambda_{1}
+b^2-2b^2(b^2-t^2)\lambda_{0}}, \ \ \ A_{1}(b^2,t)=\frac{\partial A(b^2,t)}{\partial b^2}.
\ee
Moreover, $F$ is not a Berwald metric if and only if $c_{1}\neq0$.
\end{thm}
In condition (\ref{conformal}), $\beta$ can be expressed by \cite{Z3} when $\alpha$ is of constant sectional curvature.
One can give the explicit expression of $\phi$ by choosing suitable functions $c_{1}(b^2)$, $\lambda_{0}(b^2)$, $\lambda_{1}(b^2)$ and $\lambda_{3}(b^2)$ . When $c_{1}(b^2)$, $\lambda_{0}(b^2)$, $\lambda_{1}(b^2)$, $\lambda_{3}(b^2)$ are constants, then $\phi$ is just $(4)$ in \cite{Z1}.
Thus, some explicit examples can be constructed. In Section 6, we give a new almost regular Landsberg metric.

By Theorem \ref{thm1.2}, when $c_{1}=0$, $F$ is a Berwald metric. In fact,
\be
\phi=\lambda_{3}e^{\int^s_{0}{\frac{b^2t\lambda_{1}+2b^2t\lambda_{0}}{b^2t^2\lambda_{1}+b^2-2b^2(b^2-t^2)\lambda_{0}}}dt}=
\lambda_{3}\sqrt{s^2\lambda_{1}+1-2(b^2-s^2)\lambda_{0}}.
\ee
Then $F=\alpha\phi$ is a Riemannian metric.
\begin{cor}\label{cor1.3}
Let $F=\alpha \phi(b^2,\frac{\beta}{\alpha})$ be a regular general $(\alpha,\beta)$-metric on an n-dimensional manifold with $n\geq3$. Suppose that $\beta$ is closed and conformal to $\alpha$, i.e.,
\be b_{i|j}=c a_{ij},\label{bij2}
\ee
where $c=c(x)\neq0$ is a scalar function. Then $F$ is a Landsberg metric if and only if $F$ is a Riemannian metric.
\end{cor}

Then it is natural to consider the general case when $\beta$ is not conformal to $\alpha$. It would be studied in our next paper.

\section{Berwald curvature and Landsberg curvature}
For a Finsler metric $F=F(x,y)$ on a manifold $M$, the spray $G = y^i \frac{\partial}{\partial x^i} - 2 G^i \frac{\partial}{\partial y^i}$ is a vector field on $TM$, where $G^i=G^i(x,y)$ are defined by
\be
G^i = \frac{1}{4} g^{il}\Big{\{} [F^2]_{x^m y^l} y^m -[F^2]_{x^l} \Big{\}}, \label{1}
\ee
where $g_{ij}:=\frac{1}{2}[F^2]_{y^i y^j}$ and $(g^{ij}):=(g_{ij})^{-1}$.

The following lemma is to ensure the positivity of the general $(\alpha,\beta)$-metric  \cite{Y1}.
\begin{lem}\label{lem2.1} {\rm (\cite{Y1})} Let $M$ be an n-dimensional manifold. $F=\alpha \phi(b^2,\frac{\beta}{\alpha})$ is a Finsler metric on $M$ for any Riemannian metric $\alpha$ and 1-form $\beta$ with $\|\beta\|_\alpha < b_0$ if and only if $\phi=\phi(b^2,s)$ is a positive $C^{\infty}$ function satisfying
\be \phi-s \phi_2 >0, \ \  \phi-s \phi_2 + (b^2-s^2)\phi_{22}>0,\label{13}
\ee
when $n\geq3$ or $$\phi-s \phi_2 + (b^2-s^2)\phi_{22}>0,$$
when $n=2$, where $s$ and $b$ are arbitrary numbers with $|s|\leq b <b_0.$
\end{lem} \ \

In this paper, we denote the covariant derivative of the 1-form $\beta$ with respect to the  Riemannian metric $\alpha$ by $b_{i|j}$. Moreover, for simplicity, let
$$r_{ij}=\frac{1}{2}(b_{i|j}+b_{j|i}), \ s_{ij}=\frac{1}{2}(b_{i|j}-b_{j|i}), \  r_{00}= r_{ij}y^iy^j, \ s^i_{ \ 0}=a^{ij}s_{jk}y^k,$$
$$r_i= b^j r_{ji}, \ s_i=b^js_{ji}, \ r_0 = r_iy^i, \ s_0 =s_iy^i, \ r^i=a^{ij}r_j, \ s^i=a^{ij}s_j, \ r=b^ir_i.$$

The Berwald tensor $B=B_{jkl}^i\partial_{i}\otimes dx^j \otimes dx^k \otimes dx^l$ is defined by
\be
B_{jkl}^i:=\frac{\partial^3G^i}{\partial y^j \partial y^k \partial y^l}.\label{Bjkl}
\ee
A Finsler metric is called a \emph{Berwald metric} if $B_{jkl}^i=0$, i.e. the spray coefficients $G^i=G^i(x,y)$ are quadratic in $y\in T_{x}M$ at every point $x\in M$.

 The Landsberg tensor $L = L_{ijk} dx^i \otimes dx^j \otimes dx^k$ is defined by
\be
L_{jkl}:= -\frac{1}{2} FF_{y^i} [G^i]_{y^j y^k y^l}. \label{Ljkl}
\ee
A Finsler metric is called a \emph{Landsberg metric} if $L_{jkl}=0$. Clearly, any Berwald metric is a Landsberg metric.
The {\it mean Landsberg curvature}
${\bf J}= J_k dx^k$ is defined by
\[ J_k = g^{ij} L_{ijk}.\]
A Finsler metric is called a \emph{weak Landsberg metric} if $L_{jkl}=0$.
Some results about weak landsberg metric can be found in \cite{X1, LS2}.

The geodesic coefficients $G^i$ of a general $(\alpha,\beta)$-metric $F=\alpha \phi(b^2,\frac{\beta}{\alpha})$ were given in \cite{Y1} as the following
\be
G^i = G_{\alpha}^i + Py^i + Q^i, \label{GPPP}
\ee
where
$$P = \Big\{\Theta(-2 \alpha Q s_0+r_{00}+2 \alpha^2 R r)+\alpha\Omega(r_0+s_0)\Big\}{\alpha^{-1}},$$
$$Q^i = \alpha Q s^i_{\ 0}-\alpha^2 R (r^i+s^i)+\Big\{\Psi(-2\alpha Q s_0+r_{00}+2\alpha^2 R r)+\alpha \Pi(r_0+s_0)\Big\}b^i,$$
where
$$Q=\frac{\phi_{2}}{\phi-s \phi_{2}},\ \ \Theta=\frac{(\phi-s \phi_{2})\phi_{2}-s \phi \phi_{22}}{2 \phi [\phi-s \phi_{2}+(b^2-s^2)\phi_{22}]},\ \ \Psi=\frac{\phi_{22}}{2 [\phi-s \phi_{2}+(b^2-s^2)\phi_{22}]},$$
$$R=\frac{\phi_{1}}{\phi-s \phi_{2}}, \ \ \Pi=\frac{(\phi-s \phi_{2})\phi_{12}-s \phi_{1} \phi_{22}}{(\phi-s \phi_{2}) [\phi-s \phi_{2}+(b^2-s^2)\phi_{22}]},\ \ \Omega=\frac{2 \phi_{1}}{\phi}-\frac{s \phi+(b^2-s^2)\phi_{2}}{\phi}\Pi.$$
Obviously, if $\beta$ is parallel with respect to $\alpha$ ($r_{ij}=0$ and $s_{ij}=0$), then $P=0$ and $Q^i=0$. In this case, $G^i=G^{i}_{\alpha}$ are quadratic in $y$, and $F$ is a Berwald metric.\\

In the following propositions, we  give the Berwald curvature of general $(\alpha,\beta)$-metrics and the Landsberg curvature of general $(\alpha,\beta)$-metrics. In \cite{H1}, H. Zhu obtained the expression of $B_{jkl}^i$. Here we give a different version expressed by $h_{jk}$, $h_{j}$ and our $X$ is $E$ in \cite{H1}. Here $X_1 = \frac{\partial X}{\partial b^2} $, $X_2 = \frac{\partial X}{\partial s} $, $H_1 = \frac{\partial H}{\partial b^2} $, $H_2 = \frac{\partial H}{\partial s} $ and so on.

\begin{prop} Let $F=\alpha \phi(b^2,s)$ be a general $(\alpha,\beta)$-metric on an n-dimensional manifold with $n\geq2$. Suppose that $\beta$ satisfies (\ref{conformal}), then the  Berwald curvature of $F$ is given by
\begin{equation}
\begin{split}
B_{jkl}^i  =  & \frac{1}{\alpha^4}\Big\{(H_{2}-sH_{22})[h_{jk}h_{l}+h_{jl}h_{k}+h_{kl}h_{j}]+H_{222}h_{j}h_{k}h_{l}\Big\}b^i\\&
-\frac{1}{\alpha^5}\Big\{(X-sX_{2})[h_{jk}h_{l}+h_{jl}h_{k}+h_{kl}h_{j}]+X_{22}[h_{j}h_{k}y_{l}+h_{j}h_{l}y_{k}+h_{k}h_{l}y_{j}
]\\&+sX_{22}[h_{jk}h_{l}+h_{jl}h_{k}+h_{kl}h_{j}]-X_{222}h_{j}h_{k}h_{l}\Big\}y^i+\frac{1}{\alpha^3}[(X-sX_{2})h_{kl}+
X_{22}h_{k}h_{l}]\delta_{j}^i\\&+\frac{1}{\alpha^3}[(X-sX_{2})h_{jl}+
X_{22}h_{j}h_{l}]\delta_{k}^i+\frac{1}{\alpha^3}[(X-sX_{2})h_{jk}+
X_{22}h_{j}h_{k}]\delta_{l}^i,\label{Gjkl}
\end{split}
\end{equation}
where

$$h_{j} = \alpha b_{j} - s y_{j}, \ \ h_{jk} =  \alpha^2 a_{jk} - y_{j}y_{k},$$
\be X=\frac{\phi_{2}+2s\phi_{1}}{2\phi}-H\frac{s\phi+(b^2-s^2)\phi_{2}}{\phi},\label{0X}
\ee
\be H=\frac{\phi_{22}-2(\phi_{1}-s\phi_{12})}{2[\phi-s\phi_{2}+(b^2-s^2)\phi_{22}]}.\label{0H}
\ee
\end{prop}
{\bf Proof}: By (\ref{conformal}), we get
\be r_{00}=c\alpha^2,\ \ r_{0}=c\beta,\ \ r=cb^2,\ \ r^i=cb^i,\ \ s^i=0, \ \ s^i_{\ 0}=0, \ \ s_{0}=0.\label{Dbij}
\ee
Substituting (\ref{Dbij}) into (\ref{GPPP}) yields
\begin{equation}
\begin{split}
G^i & = G_{\alpha}^i + c\alpha\{\Theta(1+2b^2R)+s\Omega\}y^i+c\alpha^2\{\Psi(1+2b^2R)+s\Pi-R\}b^i
\\&=G_{\alpha}^i+c\alpha X y^i+c\alpha^2 Hb^i.
\end{split}
\end{equation}
where $X$ and $H$ are given by (\ref{0X}) and (\ref{0H}) respectively.
By (\ref{Bjkl}) and a direct computation, we obtain (\ref{Gjkl}).\qed \\

Based on (\ref{Gjkl}) we can give the Landsberg curvature.
\begin{prop} Let $F=\alpha \phi(b^2,s)$ be a general $(\alpha,\beta)$-metric on an n-dimensional manifold with $n\geq2$. Suppose that $\beta$ satisfies (\ref{conformal}), then the  Landsberg curvature of $F$ is given by
\be
L_{jkl}=-\frac{\rho}{6 \alpha^5} \Big\{h_{j}h_{k}C_{l} + h_{j}h_{l}C_{k} + h_{k}h_{l}C_{j} + 3E_{j}h_{kl} + 3E_{k}h_{jl} + 3E_{l}h_{jk}\Big\},\label{ZYLjkl}
\ee
where

$$h_{j} = \alpha b_{j} - s y_{j}, \ \ h_{jk} =  \alpha^2 a_{jk} - y_{j}y_{k},\ \  \rho  =  \phi(\phi-s\phi_{2}),$$
\be C_{j}  =  c\alpha^2\Big\{[b^2Q+s]H_{222}+[1+sQ]X_{222}+3QX_{22}\Big\}h_{j},\label{Cj}
\ee
\be
E_{j}  =  c\alpha^2\Big\{[b^2Q+s](H_{2}-sH_{22})-[1+sQ]sX_{22}+Q[X-sX_{2}]\Big\}h_{j},\label{Ej}
\ee
\be X=\frac{\phi_{2}+2s\phi_{1}}{2\phi}-H\frac{s\phi+(b^2-s^2)\phi_{2}}{\phi},\label{X}
\ee
\be H=\frac{\phi_{22}-2(\phi_{1}-s\phi_{12})}{2[\phi-s\phi_{2}+(b^2-s^2)\phi_{22}]}.\label{H}
\ee
\end{prop}
{\bf Proof}: Substituting $F=\alpha\phi(b^2,s)$ into
 (\ref{Ljkl}) yields
\be
L_{jkl}=-\frac{1}{2} \alpha \phi [\alpha\phi]_{y^i} [G^i]_{y^j y^k y^l}=-\frac{1}{2}\rho(\alpha  Qb_{i}+ y_{i})B_{jkl}^i,\label{LjklYS}
\ee
where $\rho=\phi(\phi-s\phi_{2})$.
By plugging (\ref{Gjkl}) into (\ref{LjklYS}) and a direct computation,  we obtain (\ref{ZYLjkl}).\qed \\

\section{Berwald $(\alpha,\beta)$-metric}
In this section, we are going to give the equivalent equations of Berwald general $(\alpha,\beta)$-metrics first. Then
 Theorem \ref{thm1.3} can be proved.

\begin{lem}\label{lem6}
Let $F=\alpha \phi(b^2,s)$ be a non-Riemannian general $(\alpha,\beta)$-metric on an n-dimensional manifold $n\geq2$.
Suppose $J=J(b^2,s)$ and $M=M(b^2,s)$ are arbitrary $C^{\infty}$ functions, then the following facts hold:
\begin{itemize}
\item[(i)] $(n\geq3)$ $h_{j}h_{k}h_{l}J+h_{jk}h_{l}M+h_{jl}h_{k}M+h_{kl}h_{j}M=0$ if and only if $J=0$ and $M=0$;
\item[(ii)] $(n=2)$ $h_{j}h_{k}h_{l}J+h_{jk}h_{l}M+h_{jl}h_{k}M+h_{kl}h_{j}M=0$ if and only if $(b^2-s^2)J+3M=0$.
\end{itemize}
\end{lem}
{\bf Proof}: {\bf (i)$(n\geq3)$  }

{\bf "Necessity"}:
Contracting $h_{j}h_{k}h_{l}J+h_{jk}h_{l}M+h_{jl}h_{k}M+h_{kl}h_{j}M=0$ with $b^j$, $b^k$ and $b^l$ yields
\be (b^2-s^2)J+3M=0.\label{JM}
\ee
Let $\omega_{jk}:=(b^2-s^2)h_{jk}-h_{j}h_{k}$. Then (\ref{JM}) is equivalent to the following equation:
\be
\omega_{jk}h_{l}M+\omega_{jl}h_{k}M+\omega_{kl}h_{j}M=0.\label{MMM}
\ee
Contracting (\ref{MMM}) with $b^l$ yields
\be
\alpha(b^2-s^2)\omega_{jk}M=0.\label{bsjkm}
\ee
Contracting (\ref{bsjkm}) with $a^{jk}$ yields
\be
a^{jk}\alpha(b^2-s^2)\omega_{jk}M=(n-2)\alpha^3(b^2-s^2)^2M=0.
\ee
Since $n\geq3$, we obtain $M=0$. Then $J=0$ by (\ref{JM}).
.\\
{\bf "Sufficiency"}: It is obvious.\\

{\bf (ii) $(n=2)$.  }

{\bf "Necessity"}:
  Contracting $h_{j}h_{k}h_{l}J+h_{jk}h_{l}M+h_{jl}h_{k}M+h_{kl}h_{j}M=0$ with $b^j$, $b^k$ and $b^l$ yields
\be (b^2-s^2)J+3M=0.
\ee
{\bf "Sufficiency"}: By $3M=-(b^2-s^2)J$, we get
\begin{eqnarray*}
&h_{j}h_{k}h_{l}J+h_{jk}h_{l}M+h_{jl}h_{k}M+h_{kl}h_{j}M=\\&
\frac{1}{3}\Big\{[h_{j}h_{k}-(b^2-s^2)h_{jk}]h_{l}+[h_{k}h_{l}
-(b^2-s^2)h_{kl}]h_{j}+[h_{j}h_{l}-(b^2-s^2)h_{jl}]h_{k}\Big\}J.
\end{eqnarray*}
When $n=2$, $(b^2-s^2)h_{jk}-h_{j}h_{k}=0$. Then we get
$h_{j}h_{k}h_{l}J+h_{jk}h_{l}M+h_{jl}h_{k}M+h_{kl}h_{j}M=0$.\\
\qed \\

Based on Lemma \ref{lem6}, we can prove the following proposition.\\

\begin{prop}\label{prop_Berwald}
Let $F=\alpha \phi(b^2,\frac{\beta}{\alpha})$ be a  general $(\alpha,\beta)$-metric on an n-dimensional manifold with $n\geq3$. Suppose that $\beta$ satisfies
\[ b_{i|j}=c a_{ij},
\]
where $c=c(x)\neq0$ is a scalar function.
Then $F$ is a Berwald metric if and only if
  $\phi=\phi(b^2,\frac{\beta}{\alpha})$ satisfies
\be
H_{2}-sH_{22}=0, \ \ \ X-sX_{2}=0.\label{H2sH22XsX2}
\ee

\end{prop}

{\bf Proof :}\ \ \
Contracting (\ref{Gjkl}) with $b_{i}$ and contracting (\ref{Gjkl}) with $y_{i}$ respectively yields
\be
\frac{1}{\alpha^4}\Big\{h_{j}h_{k}h_{l}\bar{J}+h_{jk}h_{l}\bar{M}+h_{jl}h_{k}\bar{M}+h_{kl}h_{j}\bar{M}\Big\}= B_{jkl}^i b_{i}\label{JM1}
\ee
and
\be
\frac{1}{\alpha^3}\Big\{h_{j}h_{k}h_{l}\tilde{J}+h_{jk}h_{l}\tilde{M}+h_{jl}h_{k}\tilde{M}+h_{kl}h_{j}\tilde{M}\Big\}=  B_{jkl}^i y_{i},\label{JM2}
\ee
where
\be \bar{J}=b^2H_{222}+sX_{222}+3X_{22},
\ \ \ \bar{M}=X-sX_{2}-s^2X_{22}+b^2(H_{2}-sH_{22}).
\ee
\be \tilde{J}=X_{222}+sH_{222},
\ \ \ \tilde{M}=s(H_{2}-sH_{22})-sX_{22}, \label{M2}
\ee

{\bf "Necessity"}:
 By Lemma \ref{lem6},  $B_{jkl}^i=0$  if and only if $\bar{J}=\bar{M}=\tilde{J}=\tilde{M}=0$.
Then by (\ref{M2}),
\be
H_{2}-sH_{22}=X_{22}.\label{H2H22X22}
\ee
$ s \bar{M} - b^2 \tilde{M} =0$  is equivalent to
\be
X-sX_{2}=-(b^2-s^2)X_{22}.\label{XX2X22}
\ee
$\bar{J}- s \tilde{J}  =0$  is equivalent to
\be
(b^2-s^2)H_{222}=-3X_{22}.\label{H222X22}
\ee
Contracting $B_{jkl}^i$ for $i$ and $j$, we get
\begin{equation}
\begin{split}
B_{mkl}^m  =  & \frac{1}{\alpha^3}\Big\{(H_{2}-sH_{22})[2h_{k}h_{l}+(b^2-s^2)h_{kl}]+(b^2-s^2)H_{222}h_{k}h_{l}
\\&-X_{22}h_{k}h_{l}+(n+2)[(X-sX_{2})h_{kl}+X_{22}h_{k}h_{l})]\Big\}=0,\label{Gjkl6}
\end{split}
\end{equation}
Plugging (\ref{H2H22X22}), (\ref{XX2X22}) and (\ref{H222X22}) into (\ref{Gjkl6}) yields
\begin{equation}
\begin{split}
B_{mkl}^m & =   \frac{1}{\alpha^3}\Big\{X_{22}[2h_{k}h_{l}+(b^2-s^2)h_{kl}]-3X_{22}h_{k}h_{l}
-X_{22}h_{k}h_{l}+(n+2)[-(b^2-s^2)X_{22}h_{kl}+X_{22}h_{k}h_{l})]\Big\}\\&=
\frac{1}{\alpha^3}X_{22}\Big\{nh_{k}h_{l}-(n+1)(b^2-s^2)h_{kl}\Big\}=0,\label{Gjkl66}
\end{split}
\end{equation}
Contracting (\ref{Gjkl66}) with $a^{kl}$ yields
\be
\frac{1}{\alpha}(b^2-s^2)(n-n^2+1)X_{22}=0.\label{nX22}
\ee
Because $n-n^2+1\neq0$, it follows from (\ref{nX22}) that
\be
X_{22}=0.\label{X220}
\ee
Plugging (\ref{X220}) into (\ref{H2H22X22}) and (\ref{XX2X22}) yields
\be
H_{2}-sH_{22}=0, \ \ \   X-sX_{2}=0.
\ee

{\bf "Sufficiency"}:
By (\ref{H2sH22XsX2}), we get
\be
X_{22}=0, \ \ \ X_{222}=0, \ \ \ H_{222}=0.\label{XXH}
\ee

Plugging (\ref{H2sH22XsX2}) and (\ref{XXH}) into (\ref{Gjkl}) yields $B_{jkl}^i=0$. \qed

{\bf Proof of Theorem {\ref{thm1.3}}:}\ \ \
In \cite{H1}, H. Zhu gave the general solutions of (\ref{H2sH22XsX2}). Some of them are just Riemannian metrics and the others
 can be expressed by
\[ \phi=\varphi(\frac{s^2}{e^{\int{(\frac{1}{b^2}-b^2t )}db^2}+s^2\int{t e^{\int{(\frac{1}{b^2}-b^2t )}db^2}db^2}})
e^{\int{(\frac{1}{2}b^2t -\frac{1}{b^2})}db^2} s,\]
where $\varphi$ is any positive differentiable function and $t$ is a smooth function of $b^2$.
It is easy to see that $\phi$ is an odd function in $s$. When $s=0$, $\phi=0$. Thus it is excluded. Thus there is no
non-Riemannian Berwald general $(\alpha,\beta)$-metric when  the dimension $n\geq3$. \qed

\section{Landsberg $(\alpha,\beta)$-metric}
In this section, we are going to give the sufficient and necessary conditions of a general $(\alpha,\beta)$-metric to be a Landsberg metric.
By definition, a Finsler metric is called a Landsberg metric if $L_{jkl}=0$. First we  prove a equivalent condition of $L_{jkl}=0$.
Note that
$$h_{j}b^j=\alpha(b^2-s^2),\ \ \ h_{j}y^j=0,\ \ \ C_{j}y^j=0,\ \ \ E_{j}y^j=0,$$
$$h_{jk}b^k=\alpha h_{j},\ \ \ h_{jk}y^k=0,\ \ \ h_{jk}b^jb^k=\alpha^2(b^2-s^2).$$
Let
\begin{eqnarray*}
C: & = & b^j C_{j} = c\alpha^3(b^2-s^2)\Big\{[b^2Q+s]H_{222}+[1+sQ]X_{222}+3QX_{22}\Big\},\\
E: & = & b^j E_{j} = c\alpha^3(b^2-s^2)\Big\{[b^2Q+s](H_{2}-sH_{22})-[1+sQ]sX_{22}+Q[X-sX_{2}]\Big\}.
\end{eqnarray*}

\begin{lem}\label{lem2.2}
Let $F=\alpha \phi(b^2,s)$ be a general $(\alpha,\beta)$-metric on an n-dimensional manifold $n\geq2$.
Suppose that $\beta$ satisfies (\ref{conformal}), then the following hold:
\begin{itemize}
\item[(i)] $(n\geq3)$ $L_{jkl}=0$  if and only if $C_{j}=0$ and $E_{j}=0$;
\item[(ii)] $(n=2)$ $L_{jkl}=0$   if and only if $(b^2-s^2)C_{j}+3E_{j}=0$.
\end{itemize}
\end{lem}
{\bf Proof}: {\bf (i) $(n\geq3)$}.

{\bf  "Necessity"}:
Assume that $L_{jkl}=0$. Then by (\ref{ZYLjkl}), we get
\be
h_{j}h_{k}C_{l} + h_{j}h_{l}C_{k} + h_{k}h_{l}C_{j} + 3E_{j}h_{kl} + 3E_{k}h_{jl} + 3E_{l}h_{jk}=0.\label{Ljkl1}
\ee
Contracting (\ref{Ljkl1}) with $b^j$, $b^k$ and $b^l$ yields
\[
(b^2-s^2)C+3E=0.
\]
Contracting (\ref{Ljkl1}) with $b^k$ and $b^l$ yields
\be
(b^2-s^2)C_{j}+3E_{j}=0.\label{CE}
\ee
Let $\omega_{jk}:=(b^2-s^2)h_{jk}-h_{j}h_{k}$. By (\ref{CE}), (\ref{Ljkl1}) is equivalent to the following equation:
\be
\omega_{jk}C_{l}+\omega_{jl}C_{k}+\omega_{kl}C_{j}=0.\label{C}
\ee
Contracting (\ref{C}) with $b^l$ yields
\be
\omega_{jk}C=0.\label{wjkc}
\ee
Noting $a^{kl}\omega_{kl}=(n-2)\alpha^2(b^2-s^2)$ and contracting (\ref{wjkc}) with $a^{jk}$ yields
\[
a^{jk}\omega_{jk}C=(n-2)\alpha^2(b^2-s^2)C=0.
\]
Since $n\geq3$, we obtain  $C=0$. Contracting (\ref{C}) with $a^{kl}$   yields
\[
- 2 \alpha C h_{j} +  (n-1)\alpha^2(b^2-s^2)C_{j} = 0.
\]
Then $C_{j}=0$ and $E_{j}=0$ by (\ref{CE}).\\
{\bf "Sufficiency"}. It is obvious by (\ref{ZYLjkl}).\\
{\bf (ii) $(n=2)$}.

{\bf "Necessity"}:
Assume that $L_{jkl}=0$.  Then by (\ref{ZYLjkl}), we get
\be
h_{j}h_{k}C_{l} + h_{j}h_{l}C_{k} + h_{k}h_{l}C_{j} + 3E_{j}h_{kl} + 3E_{k}h_{jl} + 3E_{l}h_{jk}=0.\label{Ljkl12}
\ee
Contracting (\ref{Ljkl12}) with $b^j$, $b^k$ and $b^l$ yields
\[
(b^2-s^2)C+3E=0.
\]
Contracting (\ref{Ljkl12}) with $b^k$ and $b^l$ yields
\[
(b^2-s^2)C_{j}+3E_{j}=0.
\]
\\
{\bf "Sufficiency"}: Plugging $3E_{j}=-(b^2-s^2)C_{j}$ into (\ref{ZYLjkl}), we get
\begin{equation}
\begin{split}
L_{jkl}&=-\frac{\rho}{6 \alpha^5} \Big\{h_{j}h_{k}C_{l} + h_{j}h_{l}C_{k} + h_{k}h_{l}C_{j} -(b^2-s^2)C_{j}h_{kl}  -(b^2-s^2)C_{k}h_{jl} -(b^2-s^2)C_{l}h_{jk}\Big\}\\&
=-\frac{\rho}{6 \alpha^5} \Big\{[h_{j}h_{k}-(b^2-s^2)h_{jk}]C_{l} + [h_{j}h_{l}-(b^2-s^2)h_{jl}]C_{k} + [h_{k}h_{l}-(b^2-s^2)h_{kl}]C_{j}\Big\}.\label{Ljklzm}
\end{split}
\end{equation}
When $n=2$, $(b^2-s^2)h_{jk}-h_{j}h_{k}=0$. Then by (\ref{Ljklzm}), we get
$L_{jkl}=0$.
\qed \\

The following lemma is needed in Proposition \ref{thm1.1}.
\begin{lem}\label{SUF}
Let $F=\alpha \phi(b^2,s)$ be a  general $(\alpha,\beta)$-metric on an n-dimensional manifold $n\geq 2$. Suppose that $\beta$ satisfies
\be
b_{i|j}=c a_{ij},\label{bij0}
\ee
where $c=c(x)\neq0$ is a scalar function.
If one of the following holds:\\
(i) $(n\geq3)$  The function $\phi=\phi(b^2,\frac{\beta}{\alpha})$ satisfies
\be
[b^2Q+s](H_{2}-sH_{22})-[1+sQ]sX_{22}+Q[X-sX_{2}]=0\label{ZYFC10}
\ee
and
\be
[b^2Q+s]H_{222}+[1+sQ]X_{222}+3QX_{22}=0.\label{ZYFC20}
\ee
(ii) $(n=2)$ The function $\phi=\phi(b^2,\frac{\beta}{\alpha})$ satisfies
\be
(b^2-s^2)\Big\{[b^2Q+s]H_{222}+[1+sQ]X_{222}+3QX_{22}\Big\}
+3\Big\{[b^2Q+s](H_{2}-sH_{22})-[1+sQ]sX_{22}+Q[X-sX_{2}]\Big\}=0.\label{bsC3E}
\ee

then $F=\alpha \phi(b^2,s)$ is a Landsberg metric.
\end{lem}
{\bf Proof}:
When the dimension $n\geq3$, by the assumption  (\ref{ZYFC10}) and (\ref{ZYFC20}), we get
$$
C_{j} = c\alpha^2\Big\{[b^2Q+s]H_{222}+[1+sQ]X_{222}+3QX_{22}\Big\}h_{j}=0,
$$
$$
E_{j} = c\alpha^2\Big\{[b^2Q+s](H_{2}-sH_{22})-[1+sQ]sX_{22}+Q[X-sX_{2}]\Big\}h_{j}=0.
$$
By Lemma \ref{lem2.2}, we obtain $L_{jkl}=0$.

When the dimension $n=2$, by $(\ref{bsC3E})\times c\alpha^2h_{j}$, we get
\begin{equation}
\begin{split}
&(b^2-s^2)c\alpha^2\Big\{[b^2Q+s]H_{222}+[1+sQ]X_{222}+3QX_{22}\Big\}h_{j}\\&
+3c\alpha^2\Big\{[b^2Q+s](H_{2}-sH_{22})-[1+sQ]sX_{22}+Q[X-sX_{2}]\Big\}h_{j}=0.
\end{split}
\end{equation}
Which is equivalent to
\be
(b^2-s^2)C_{j}+3E_{j}=0.
\ee
By Lemma \ref{lem2.2}, we obtain $L_{jkl}=0$.
\qed

Now we can prove the following proposition.

\begin{prop}\label{thm1.1}
Let $F=\alpha \phi(b^2,\frac{\beta}{\alpha})$ be a  general $(\alpha,\beta)$-metric on an n-dimensional manifold with $n\geq2$. Suppose that $\beta$ satisfies
\be b_{i|j}=c a_{ij},\label{bij}
\ee
where $c=c(x)\neq0$ is a scalar function.
Then $F$ is a Landsberg metric if and only if one of the following holds:
\\
(i) When the dimension $n\geq3$, the function $\phi=\phi(b^2,\frac{\beta}{\alpha})$ satisfies
\be
[b^2Q+s](H_{2}-sH_{22})-[1+sQ]sX_{22}+Q[X-sX_{2}]=0,\label{ZYFC1}
\ee
\be
[b^2Q+s]H_{222}+[1+sQ]X_{222}+3QX_{22}=0.\label{ZYFC2}
\ee
(ii)  When the dimension $n=2$, the function $\phi=\phi(b^2,\frac{\beta}{\alpha})$ satisfies
\be \label{ZYFC3}
(b^2-s^2)\Big\{[b^2Q+s]H_{222}+[1+sQ]X_{222}+3QX_{22}\Big\}
+3\Big\{[b^2Q+s](H_{2}-sH_{22})-[1+sQ]sX_{22}+Q[X-sX_{2}]\Big\}=0.
\ee

\end{prop}

{\bf Proof :}  The sufficiency is proved by Lemma \ref{SUF} . We only need to prove the necessity. When $n\geq3$, by Lemma \ref{lem2.2}, we get $C_{j}=0$ and $E_{j}=0$.
Then $\phi=\phi(b^2,s)$ satisfies (\ref{ZYFC1}) and (\ref{ZYFC2}) by (\ref{Cj}) and (\ref{Ej}).
When $n=2$, by Lemma \ref{lem2.2}, we get $(b^2-s^2)C_{j}+3E_{j}=0$. Then by (\ref{Cj}) and (\ref{Ej}), we get
(\ref{ZYFC3}).
\qed

\section{Proof of Theorem \ref{thm1.2}}
In this section, we give the general solutions of (\ref{ZYFC1}) and (\ref{ZYFC2}) to prove Theorem \ref{thm1.2}.
The following lemmas are needed in our proof.
\begin{lem}\label{lem2}
If a $C^{\infty}$ function $\phi=\phi(b^2,s)$ satisfies $Q=0$ (i.e. $\frac{\phi_{2}}{\phi-s\phi_{2}}=0$), then
\[  \phi=f(b^2),
\]
where $f(b^2)$ is an arbitrary $C^{\infty}$ function of $b^2$.
\end{lem}
{\bf Proof}: By the assumption,
\[ Q=\frac{\phi_{2}}{\phi-s\phi_{2}}=0.
\]
Then $\phi_{2}=0$.\qed

\begin{lem}\label{lem1}
If a positive $C^{\infty}$ function $\phi=\phi(b^2,s)$ satisfies $Q-sQ_{2}=0$, then
\be
\phi(b^2,s)=g(b^2)\sqrt{1+h(b^2)s^2},\label{fc1}
\ee
where $g(b^2)>0$, $h(b^2)$ are $C^{\infty}$ functions of $b^2$.
\end{lem}
{\bf Proof}:
Differentiating $Q-sQ_{2}=0$ with respect to $s$ yields
$-sQ_{22}=0$. Then we get
\be Q_{22}=0.\label{Q220}
\ee
Integrating (\ref{Q220}) with respect to $s$ yields
\be Q=h(b^2)s+f(b^2),\label{Qhf}
\ee
where $h(b^2)$ and $f(b^2)$ are $C^{\infty}$ functions of $b^2$.
Substitute (\ref{Qhf}) into $Q-sQ_{2}=0$ yields $f(b^2)=0$. Then
\be Q=h(b^2)s.\label{Qh}
\ee
By $Q=\frac{\phi_{2}}{\phi-s\phi_{2}}$, we get
\be
\frac{\phi_{2}}{\phi}=\frac{Q}{1+sQ}.
\ee
Which is equivalent to
\be [\ln\phi]_{2}=\frac{Q}{1+sQ}.
\ee
Integrating the above equation with respect to $s$ yields
\be
\phi=g(b^2)e^{\int{\frac{Q}{1+sQ}}ds}.\label{J}
\ee
Substitute (\ref{Qh}) into (\ref{J}) yields (\ref{fc1}).
\qed \\

\begin{lem}\label{lem4.3}
Let $\phi=\phi(b^2,s)$ is a positive $C^{\infty}$ function. Then $\phi$ satisfies (\ref{ZYFC1}) and (\ref{ZYFC2}) if and only if the function $\phi=\phi(b^2,s)$ satisfies
\be X-sX_{2}=\frac{c_{1}}{\sqrt{b^2-s^2}},\label{XX}
\ee
\be H_{2}-sH_{22}=-\frac{c_{1}}{(b^2-s^2)^{\frac{3}{2}}},\label{HH}
\ee
where $c_{1}=c_{1}(b^2)$ is a $C^{\infty}$ functions of $b^2$.
\end{lem}
{\bf Proof:} {\bf "Necessity"}: Differentiating (\ref{ZYFC1}) with respect to $s$ yields
\begin{equation}
\begin{split}
&[b^2Q+s](-sH_{222})+[b^2Q_{2}+1](H_{2}-sH_{22})-[Q+sQ_{2}]sX_{22}
\\&-[1+sQ]X_{22}-[1+sQ]sX_{222}+Q_{2}[X-sX_{2}]+Q[-sX_{22}]=0.\label{ZYFC160}
\end{split}
\end{equation}
By $(\ref{ZYFC2})\times s+(\ref{ZYFC160})$,
\be [b^2Q_{2}+1](H_{2}-sH_{22})-[1+s^2Q_{2}]X_{22}+Q_{2}[X-sX_{2}]=0.\label{FC1}
\ee
By $(\ref{FC1})\times Q-(\ref{ZYFC1})\times Q_{2}$,
\be (Q-sQ_{2})\Big[H_{2}-sH_{22}-X_{22}\Big]=0.\label{111}
\ee
By Lemma \ref{lem1}, we get
 $Q-sQ_{2}\neq0$ because $F$ is a non-Riemannian metric by the assumption.
Hence, it follows from (\ref{111})
\be X_{22}=H_{2}-sH_{22}.\label{222}
\ee
Differentiating (\ref{222}) with respect to $s$ yields
\be X_{222}=-sH_{222}.\label{444}
\ee
Plugging (\ref{222}) and (\ref{444}) into (\ref{ZYFC2}) yields
\be
Q\{(b^2-s^2)H_{222}+3(H_{2}-sH_{22})\}=0.\label{555}
\ee
By Lemma \ref{lem2}, we get $Q\neq0$. Then it follows from (\ref{555}) that
\be
(b^2-s^2)H_{222}+3(H_{2}-sH_{22})=0.\label{HHH}
\ee
Plugging (\ref{222}) into (\ref{ZYFC1}) yields
\be
Q\{(b^2-s^2)X_{22}+(X-sX_{2})\}=0.\label{333}
\ee
By Lemma \ref{lem2}, we get $Q\neq0$. By (\ref{333}), we get
\be (b^2-s^2)X_{22}+(X-sX_{2})=0.\label{XXX}
\ee
In fact, (\ref{HHH}) and (\ref{XXX}) can be solved.
(\ref{HHH}) is equivalent to
\[
\frac{-sH_{222}}{H_{2}-sH_{22}}=\frac{3s}{b^2-s^2}.
\]
Which can be written in
\[
\Big[\ln{|H_{2}-sH_{22}|}\Big]_{2}=\Big[-\frac{3}{2}\ln{(b^2-s^2)}\Big]_{2}.
\]
Integrating the above equation with respect to $s$ yields
\be H_{2}-sH_{22}=\frac{c_{2}}{(b^2-s^2)^{\frac{3}{2}}},\label{HH1}
\ee
where
$c_{2}=c_{2}(b^2)$ is a $C^{\infty}$ function of $b^2$.\\
By a similar argument for (\ref{XXX}), we get
\be X-sX_{2}=\frac{c_{1}}{\sqrt{b^2-s^2}},\label{XX1}
\ee
where $c_{1}=c_{1}(b^2)$ is a $C^{\infty}$ function of $b^2$.\\
Plugging (\ref{HH1}) and (\ref{XX1}) into (\ref{222}) yields
\be \frac{c_{1}}{(b^2-s^2)^{\frac{3}{2}}}+\frac{c_{2}}{(b^2-s^2)^{\frac{3}{2}}}=0.
\ee
Then
\be c_{2}=-c_{1}.\label{cc}
\ee
Substitute (\ref{cc}) into (\ref{HH1}) yields
\be
H_{2}-sH_{22}=-\frac{c_{1}}{(b^2-s^2)^{\frac{3}{2}}}.\label{HH2}
\ee
If $\phi=\phi(b^2,s)$ satisfies (\ref{XX}), (\ref{HH}) and $c_{1}\neq0$, then $F=\alpha\phi(b^2,s)$ is an almost regular Landsberg metric. This metric might be singular in two directions $y\in T_{x}M$ with $\beta(x,y)=\pm b\alpha(x,y)$.\\
{\bf "Sufficiency"}. It is obvious by plugging (\ref{XX}), (\ref{HH}) back to (\ref{ZYFC1}) and (\ref{ZYFC2}). By a
direct computation, we get
\begin{equation}
\begin{split}
&[b^2Q+s](H_{2}-sH_{22})-[1+sQ]sX_{22}+Q[X-sX_{2}]\\=&-[b^2Q+s]\frac{c_{1}}{(b^2-s^2)^{\frac{3}{2}}}
+[s+s^2Q]\frac{c_{1}}{(b^2-s^2)^{\frac{3}{2}}}+Q\frac{c_{1}}{\sqrt{b^2-s^2}}=0,
\end{split}
\end{equation}
\begin{equation}
\begin{split}
&[b^2Q+s]H_{222}+[1+sQ]X_{222}+3QX_{22}\\=&[b^2Q+s]\frac{3c_{1}}{(b^2-s^2)^\frac{5}{2}}
-[1+sQ]\frac{3sc_{1}}{(b^2-s^2)^\frac{5}{2}}-3Q\frac{c_{1}}{(b^2-s^2)^\frac{3}{2}}
=0.
\end{split}
\end{equation}
\qed\\

{\bf Proof of Theorem {\ref{thm1.2}}}:\ \ \
By Proposition \ref{thm1.1}, $\phi$ satisfies (\ref{ZYFC1}) and (\ref{ZYFC2}). Then $\phi$ satisfies (\ref{XX}) and (\ref{HH}) by Lemma \ref{lem4.3}.
Thus we only need to solve (\ref{XX}) and (\ref{HH}) for $\phi$.
Substituting (\ref{XX}), (\ref{HH}) back to (\ref{XXX}), (\ref{HHH}) yields
\[
X_{22}=-\frac{c_{1}}{(b^2-s^2)^\frac{3}{2}},
\]
\[
H_{222}=\frac{3c_{1}}{(b^2-s^2)^\frac{5}{2}}.
\]
Then by integrating $X_{22}$ and $H_{222}$ with respect to $s$, we get
\be X=\frac{c_{1}\sqrt{b^2-s^2}}{b^2}+\lambda_{1}s+k_{0},\label{X12}
\ee
\be H=\frac{1}{2}\lambda_{2}s^2-\frac{c_{1}s\sqrt{b^2-s^2}}{b^4}+k_{1}s+\lambda_{0},\label{H12}
\ee
where $c_{1}=c_{1}(b^2)$, $k_{0}=k_{0}(b^2)$, $k_{1}=k_{1}(b^2)$, $\lambda_{0}=\lambda_{0}(b^2)$, $\lambda_{1}=\lambda_{1}(b^2)$ and $\lambda_{2}=\lambda_{2}(b^2)$ are $C^{\infty}$ functions of $b^2$.
Substituting (\ref{X12}) and (\ref{H12}) into (\ref{XX}) and (\ref{HH}) yields $k_{0}=0$ and $k_{1}=0$. Then we get
\be X=\frac{c_{1}\sqrt{b^2-s^2}}{b^2}+\lambda_{1}s,\label{X1}
\ee
\be H=\frac{1}{2}\lambda_{2}s^2-\frac{c_{1}s\sqrt{b^2-s^2}}{b^4}+\lambda_{0},\label{H1}
\ee
where $c_{1}=c_{1}(b^2)$, $\lambda_{0}=\lambda_{0}(b^2)$, $\lambda_{1}=\lambda_{1}(b^2)$ and $\lambda_{2}=\lambda_{2}(b^2)$ are $C^{\infty}$ functions of $b^2$.\\
By (\ref{X}) and (\ref{H}), we get
\be 2\phi X-\phi_{2}-2s\phi_{1}+2H\Big[s\phi+(b^2-s^2)\phi_{2}\Big]=0,\label{123}
\ee
\be 2H\Big[\phi-s\phi_{2}+(b^2-s^2)\phi_{22}\Big]-\phi_{22}+2(\phi_{1}-s\phi_{12})=0.\label{456}
\ee
Differentiating (\ref{123}) with respect to $s$ yields
\be 2\phi_{2}X+2\phi X_{2}-\phi_{22}-2\phi_{1}-2s\phi_{12}+2H_{2}[s\phi+(b^2-s^2)\phi_{2}]
+2H[\phi-s\phi_{2}+(b^2-s^2)\phi_{22}]=0.\label{123456}
\ee
To cancel the term involving $\phi_{12}$, we consider $(\ref{123456})-(\ref{456})$. Then we get
\be \phi_{2}X+\phi X_{2}-2\phi_{1}+H_{2}[s\phi+(b^2-s^2)\phi_{2}]
=0.\label{XH14}
\ee
$(\ref{XH14})\times s-(\ref{123})$ yields
\be
\Big\{sX+sH_{2}(b^2-s^2)+1-2H(b^2-s^2)\Big\}\phi_{2}+\Big\{sX_{2}+s^2H_{2}-2X-2sH\Big\}\phi=0.\label{666}
\ee
Plugging (\ref{X1}) and (\ref{H1}) into (\ref{666}) yields
\begin{equation}
\begin{split}
\Big\{2sc_{1}\sqrt{b^2-s^2}+b^2s^2\lambda_{1}+b^2-2b^2(b^2-s^2)\lambda_{0}\Big\}\phi_{2}
-\Big\{2c_{1}\sqrt{b^2-s^2}+b^2s\lambda_{1}+2b^2s\lambda_{0}\Big\}\phi=0.
\end{split}
\end{equation}
Which is equivalent to
\be
\Big\{2sc_{1}\sqrt{b^2-s^2}+b^2s^2\lambda_{1}+b^2-2b^2(b^2-s^2)\lambda_{0}\Big\}\Big(\ln{\phi}\Big)_{2}
=\Big\{2c_{1}\sqrt{b^2-s^2}+b^2s\lambda_{1}+2b^2s\lambda_{0}\Big\}.
\ee
Then by a direct computation, we get
\be
\phi=\lambda_{3}e^{\int^s_{0}{\frac{2c_{1}\sqrt{b^2-t^2}+b^2t\lambda_{1}+2b^2t\lambda_{0}}{2tc_{1}\sqrt{b^2-t^2}+b^2t^2\lambda_{1}
+b^2-2b^2(b^2-t^2)\lambda_{0}}}dt}.\label{TJ0}
\ee
where $c_{1}=c_{1}(b^2)$, $\lambda_{0}=\lambda_{0}(b^2)$, $\lambda_{1}=\lambda_{1}(b^2)$ and $\lambda_{3}=\lambda_{3}(b^2)>0$ are $C^{\infty}$ functions of $b^2$.
It can be seen that $\phi(b^2,s)$ might be singular in two directions $y\in T_{x}M$ with $\beta(x,y)=\pm b\alpha(x,y)$.
To prove the sufficiency, we only need substitute (\ref{TJ0}) into (\ref{123}). In fact, $c_{1}(b^2)$, $\lambda_{0}(b^2)$, $\lambda_{1}(b^2)$ and $\lambda_{3}(b^2)$ are not arbitrary. They should satisfy the following condition,
\begin{equation}
\begin{split}
\int^s_{0}{A_{1}(b^2,t)}dt=&\frac{1}{2b^4s\lambda_{3}}\Big\{\Big[b^4(b^2-s^2)(2\lambda_{0}+s^2\lambda_{2})
-2sc_{1}(b^2-s^2)^{\frac{3}{2}}-b^4\Big]\lambda_{3}A(b^2,s)\\&+\Big[b^4s(\lambda_{2}s^2+2\lambda_{1}+2\lambda_{0})
+2c_{1}(b^2-s^2)^{\frac{3}{2}}\Big]\lambda_{3}-2b^4s\lambda_{3}'\Big\},
\end{split}
\end{equation}
where
\be
A(b^2,t)=\frac{2c_{1}\sqrt{b^2-t^2}+b^2t\lambda_{1}+2b^2t\lambda_{0}}{2tc_{1}\sqrt{b^2-t^2}+b^2t^2\lambda_{1}
+b^2-2b^2(b^2-t^2)\lambda_{0}}, \ \ \ A_{1}(b^2,t)=\frac{\partial A(b^2,t)}{\partial b^2}.
\ee

Plugging (\ref{X1}) and (\ref{H1}) into $G^i$ yields
$$
G^i  =G_{\alpha}^i+c\alpha X y^i+c\alpha^2 Hb^i=G_{\alpha}^i+c\Big\{\Big[\frac{c_{1}\sqrt{b^2\alpha^2-\beta^2}}{b^2}+\lambda_{1}\beta\Big]y^i+
\Big[\frac{1}{2}\lambda_{2}\beta^2-\frac{c_{1}\beta\sqrt{b^2\alpha^2-\beta^2}}{b^4}+\lambda_{0}\alpha^2\Big]b^i\Big\}.
$$
Thus if $c_{1}\neq0$, then $G^i$ are not quadratic in $y$. Hence $F$ is not a Berwald metric.\qed

\section{Explicit Example}
In this section, we construct a new Landsberg metric by Theorem \ref{thm1.2}. In order to find explicit metric, we need to choose suitable $c_{1}(b^2)$, $\lambda_{0}(b^2)$, $\lambda_{1}(b^2)$ and $\lambda_{3}(b^2)$.\\

Differentiating (\ref{TJTJ}) with respect to $s$ yields
\be
\varsigma(b^2)s^2-b^2\sqrt{b^2-s^2}\tau(b^2)s-b^2\varsigma(b^2)=0,\label{LZ}
\ee
where
\be
\varsigma(b^2)=\Big[2(2b^2\lambda_{0}+b^4\lambda_{2}+b^2\lambda_{1}+1)\lambda_{0}-\lambda_{1}
+b^2(4\lambda_{0}'-\lambda_{2})\Big]c_{1}+2(1-2b^2\lambda_{0})c_{1}',
\ee
\be
\tau(b^2)=(2b^2\lambda_{0}-1)\lambda_{1}^2+\Big[2b^2(b^2\lambda_{2}+2\lambda_{0})\lambda_{0}
+2b^2\lambda_{0}'-b^2\lambda_{2}\Big]\lambda_{1}+2b^2(\lambda_{2}-\lambda_{1}')\lambda_{0}+4\lambda_{0}^2-\lambda_{2}
+\lambda_{1}'+2\lambda_{0}'.
\ee
Because $s$ is a  variable  in  (\ref{LZ}), we get
\be
\varsigma(b^2)=0,\ \ \ \tau(b^2)=0.\label{TJDR}
\ee
From the above equation, if we set the two functions in $\lambda_{0}$, $\lambda_{1}$, $\lambda_{2}$, $c_{1}$, then
the other functions can be determined by (\ref{TJDR}). Then we can get $\lambda_{3}$ by (\ref{TJTJ}).

\begin{ex} Put
$\lambda_{0}(b^2)=0$, $\lambda_{1}(b^2)=0$ in (\ref{TJDR}), we get $\lambda_{2}=0$, $c_{1}=k$ ($k$ is a constant). Then we can get $\lambda_{3}=m$($m$ is a constant). Let $c_{1}=1$, $\lambda_{3}=1$, then
\be
\phi=e^{\int^s_{0}{\frac{2\sqrt{b^2-t^2}}{2t\sqrt{b^2-t^2}+b^2}}dt}.\label{TJ1}
\ee
To ensure the positivity of $F$, $\phi$ should satisfy Lemma \ref{lem2.1}. In the following we prove that $\phi$ satisfies Lemma \ref{lem2.1}.
Differentiating (\ref{TJ1}) with respect to $s$ yields
\be
\phi_{2}=e^{\int^s_{0}{\frac{2\sqrt{b^2-t^2}}{2t\sqrt{b^2-t^2}+b^2}}dt}
\Big[\frac{2\sqrt{b^2-s^2}}{2s\sqrt{b^2-s^2}+b^2}\Big].\label{TJ2}
\ee
Differentiating (\ref{TJ2}) with respect to $s$ yields
\[
\phi_{22}=-e^{\int^s_{0}{\frac{2\sqrt{b^2-t^2}}{2t\sqrt{b^2-t^2}+b^2}}dt}
\Big[\frac{2b^2s}{\sqrt{b^2-s^2}(2s\sqrt{b^2-s^2}+b^2)^2}\Big].
\]
Then
\[
\phi-s\phi_{2}=e^{\int^s_{0}{\frac{2\sqrt{b^2-t^2}}{2t\sqrt{b^2-t^2}+b^2}}dt}
\Big[1-\frac{2s\sqrt{b^2-s^2}}{2s\sqrt{b^2-s^2}+b^2}\Big]>0,
\]
\[
\phi-s\phi_{2}+(b^2-s^2)\phi_{22}=e^{\int^s_{0}{\frac{2\sqrt{b^2-t^2}}{2t\sqrt{b^2-t^2}+b^2}}dt}
\Big[1-\frac{4s^2(b^2-s^2)+4b^2s\sqrt{b^2-s^2}}{4s^2(b^2-s^2)+4b^2s\sqrt{b^2-s^2}+b^4}\Big]>0.
\]
Thus $F=\alpha\phi$ is an almost regular Landsberg (non-Berwaldian) metric when $\beta$ is closed and conformal to
$\alpha$.
\end{ex}

{}

\end{document}